\documentclass[10pt]{article}

\usepackage{amsthm,amsmath,amssymb}
\usepackage{yfonts}

\usepackage{graphicx,url}
\usepackage{tikz-cd}

\usepackage[colorlinks=true,citecolor=blue,linkcolor=blue,urlcolor=blue]{hyperref}

\newcommand{\arxiv}[1]{\href{http://arxiv.org/abs/#1}{\texttt{arXiv:#1}}}

\theoremstyle{plain}
\newtheorem{theorem}{Theorem}[section]
\newtheorem{lemma}[theorem]{Lemma}
\newtheorem{corollary}[theorem]{Corollary}

\theoremstyle{definition}
\newtheorem{definition}[theorem]{Definition}
\newtheorem{example}[theorem]{Example}
\newtheorem{conjecture}[theorem]{Conjecture}

\theoremstyle{remark}

\title{\bf Some properties of the Redei-Berge function and related combinatorial Hopf algebras }

\author{Stefan Mitrovi\'c\footnote{Corresponding author\\ Address : stefan.mitrovic@matf.bg.ac.rs}\\
\small Faculty of Mathematics\\[-0.8ex]
\small University of Belgrade\\[-0.8ex]
\small Serbia\\
\small\tt stefan.mitrovic@matf.bg.ac.rs\\
\and
Tanja Stojadinovi\'c\\
\small Faculty of Mathematics\\[-0.8ex]
\small University of Belgrade\\[-0.8ex]
\small Serbia\\
\small\tt tanja.stojadinovic@matf.bg.ac.rs}

\begin{document}

\maketitle

\begin{abstract}
Stanley and Grinberg introduced a symmetric function associated to digraphs, called the Redei-Berge symmetric function. In \cite{GS} is shown that this symmetric function arises from a suitable structure of combinatorial Hopf algebra on digraphs. In this paper, we introduce two new combinatorial Hopf algebras of posets and permutations and define corresponding Redei-Berge functions for them. By using both theories, of symmetric functions and of combinatorial Hopf algebras, we prove many properties of the Redei-Berge function. 
These include some forms of deletion property, which make it similar to the chromatic symmetric function. We also find some invariants of digraphs that are detected by the Redei-Berge function. 


\end{abstract}

\section{Introduction}

Redei's and Berge's theorems are two beautiful graph-theoretical results that deserve to be  better known than they are. Berge's theorem says that if $X$ is a simple digraph and $\overline{X}$ is its complement, then the number of Hamiltonian paths of $\overline{X}$ is congruent to the same number for $X$ modulo $2$. Redei's theorem says that if $X$ is a tournament, then the number of Hamiltonian paths of $X$ is odd.

A recent paper of Grinberg and Stanley \cite{S} has revealed these results as a consequence of symmetric function theory, by constructing a symmetric function $U_X$ (the Redei-Berge symmetric function) whose $p-$basis expansion easily yields both theorems. The function $U_X$ first appeared in the 1990s, in Chow's paper \cite{C}, and later in Wiseman's \cite{W} in 2007, in connection with the path-cycle symmetric function of a digraph.
The function $U_X$ is in \cite{GS} identified as the image of the isomorphism class of digraph $X$ under a certain canonical Hopf algebra morphism, defined through the Aguiar-Bergeron-Sottile universal construction applied to a combinatorial Hopf algebra of digraphs. This reconceptualization is used to prove some properties of $U_X$, including the antipode formula. The authors in \cite{GS} also reprove Berge's theorem by combining this antipode formula with a deletion-contraction-like recurence for the principal evaluation of $U_X$. This principal specialization, $u_X(m)$, is called the Redei-Berge polynomial. 

In this paper, which serves as a continuation of \cite{S} and \cite{GS}, we further investigate properties of the Redei-Berge function of digraphs. In Section 3, we introduce two combinatorial Hopf algebras compatible with the Redei-Berge combinatorial Hopf algebra of digraphs, which enables us to define the Redei-Berge function of permutations and posets. In the sequel, in Section 4, we investigate which invariants of digraphs and posets can be obtained from the Redei-Berge function and prove some deletion properties of this function. Finally, we construct numerous new bases for the algebra of symmetric functions whose generators are Redei-Berge functions. The reader who is not familiar with the theory of combinatorial Hopf algebras and who is mainly interested in the properties of the Redei-Berge function of digraphs can focus only on Section 4 since it does not rely on the results from Section 3.

\section{Preliminaries}

The theory of combinatorial Hopf algebras is founded in \cite{ABS}.
We review some basic notions and facts. A {\it set partition} $\{V_1,\ldots,V_k\}\vdash V$ of length $k$
of a finite set $V$ is a set of disjoint nonempty subsets with $V_1\cup\ldots\cup V_k=V$. A {\it set composition} $(V_1,\ldots,V_k)\models V$ is an ordered partition. If $A\subseteq \mathbb{N}$ is finite, the unique order preserving map $\text{st}: A\rightarrow\{1, \ldots, |A|\}$ is called the \textit{standardization} of $A$.

A \textit{composition} $\alpha\models n$, where $n\in\mathbb{N}$, is a sequence $\alpha=(a_1,\ldots,a_k)$ of
positive integers with $a_1+\cdots+a_k=n$. The \textit{type of a set composition} $(V_1,\ldots,V_k)\models V$ is the composition $\mathrm{type}(V_1,\ldots,V_k)=(|V_1|,\ldots,|V_k|)\models |V|$. A \textit{partition} $\lambda\vdash n$ is a
composition $(a_1,\ldots,a_k)\models n$ such that $a_1\geq a_2\geq\cdots\geq
a_k$. There is a bijection between the set of compositions of $n$, denoted by $\mathrm{Comp}(n)$ and the power set $2^{[n-1]}$ given by
$(a_1,\ldots,a_k)\mapsto\{a_1,a_1+a_2,\ldots, a_1+\cdots+a_{k-1}\}$. We denote the inverse of this bijection by $I\mapsto\mathrm{comp}(I)$.

A {\it combinatorial Hopf algebra}, CHA
for short, $(\mathcal{H},\zeta)$  over a field $\mathbf{k}$ is a graded, connected Hopf
algebra $\mathcal{H}=\oplus_{n\geq 0}\mathcal{H}_n$ over $\mathbf{k}$ together with
a $\mathbf{k}-$algebra homomorphism
$\zeta:\mathcal{H}\rightarrow\mathbf{k}$ called the {\it character}, \cite{ABS}.

The terminal object in the category of CHAs is the
CHA of quasisymmetric functions $(QSym,\zeta_Q)$. For basics of
quasisymmetric functions see \cite{EC} and for the Hopf algebra structure on $QSym$ see \cite{ABS}. A composition $\alpha=(a_1,\ldots,a_k)\models n$ defines the monomial quasisymmetric function

\[M_\alpha=\sum_{i_1<\cdots<
i_k}x_{i_1}^{a_1}\cdots x_{i_k}^{a_k}.\] Alternatively, we write $M_I=M_{\mathrm{comp}(I)},\ I\subseteq[n-1]$. Another basis consists of fundamental quasisymmetric functions 

\begin{equation}\label{fundamental}
F_I=\sum_{\substack{1\leq i_1\leq i_2\leq\cdots\leq i_n\\
                  i_j<i_{j+1} \ \mathrm{for \ each} \ j\in I}}x_{i_1}x_{i_2}\cdots x_{i_n}, \quad I\subseteq[n-1],
\end{equation}
which are expressed in the monomial basis by
\begin{equation}\label{montofund}
F_I=\sum_{I\subseteq J}M_J, \quad I\subseteq[n-1].
\end{equation}

The unique canonical morphism
$\Psi:(\mathcal{H},\zeta)\rightarrow(QSym,\zeta_Q)$ for a CHA $(\mathcal{H}, \zeta)$ is given on
homogeneous elements with

\begin{equation}\label{canonical}
\Psi(h)=\sum_{\alpha\models n}\zeta_\alpha(h)M_\alpha, \
h\in\mathcal{H}_n,
\end{equation} where $\zeta_\alpha$ is the convolution product
\begin{equation}\label{generalcoeff}
\zeta_\alpha=\zeta_{a_1}\cdots\zeta_{a_k}:\mathcal{H}
\stackrel{\Delta^{(k-1)}}\longrightarrow\mathcal{H}^{\otimes
k}\stackrel{proj}\longrightarrow\mathcal{H}_{a_1}\otimes\cdots\otimes\mathcal{H}_{a_k}
\stackrel{\zeta^{\otimes k}}\longrightarrow\mathbf{k}.
\end{equation}

The algebra of symmetric functions $Sym$ is a subalgebra of
$QSym$ and it is the terminal object in the category of
cocommutative combinatorial Hopf algebras. This algebra consists of quasisymmetric functions that are invariant under the action of permutations on the set of variables. As a vector space, it has many natural bases, see \cite{EC}. We will be particularly interested in the \textit{power sum basis}. The $i$th \textit{power sum symmetric function} is defined as
\[p_0=1,\]
\[p_i=\sum_{j=1}^{\infty}x_j^i,\]
for $i\geq 1$. For partition $\lambda=(\lambda_1, \lambda_2, \ldots, \lambda_k)$, we define 

\[p_{\lambda}=p_{\lambda_1}p_{\lambda_2}\cdots p_{\lambda_k}.\] The power sum functions form a basis of $Sym$ only if the background field $\mathbf{k}$ has characteristic 0, which is the case we will be interested in.

A {\it digraph} $X$ is a pair $X=(V,E)$, where $V$ is a finite set and $E$ is a collection
$E\subseteq V\times V$. Elements $u\in V$ are vertices and
elements $(u,v)\in E$ are directed edges of the digraph $X$. Note that this definition allows loops as edges of a digraph, as in \cite{S}. On the other hand, the authors in \cite{GS} forbid the loops to appear. The reader should not be concerned about this inconsistency since loops do not affect the function we will consider in this paper.

A $V$-{\it listing} is a list of all vertices with no repetitions, i.e. a bijective map $\omega:[n]\rightarrow V$. We write $\Sigma_V$ for the set of all $V$-listings. For a $V$-listing $\omega=(\omega_1,\ldots,\omega_n)\in\Sigma_V$, define the $X$-{\it descent set} as

\[X\mathrm{Des}(\omega)=\{1\leq i\leq n-1\ | \ (\omega_i,\omega_{i+1})\in E\}.\] If $i\in X\mathrm{Des}(\omega)$, we say that the edge $\omega_i\omega_{i+1}$ is detected by $\omega$. If $V$ is the set $[n]=\{1,\ldots,n\}$ and $E=\{(i,j)\ | \ 1\leq j<i\leq n\}$ then $X$-descent sets are standard descent sets of permutations $\omega\in\Sigma_{[n]}$.

If $X=(V, E)$ is a digraph, its \textit{complementary digraph} is the digraph $\overline{X}=(V, (V\times V)\setminus E)$ and its \textit{opposite digraph} is $X^{op}=(V, E')$, where $E'=\{(v, u) \ \mid \ (u, v)\in E\}$. 

Grinberg and Stanley associated to a digraph $X$ a generating function for $X$-descent sets  expanded in the basis of fundamental quasisymmetric functions

\begin{equation}\label{descents}
U_X=\sum_{\omega\in\Sigma_V}F_{X\mathrm{Des}(\omega)}
\end{equation}
and named it the Redei-Berge symmetric function, see \cite{S}, \cite{RS}.

\begin{definition}
    Let $X=(V, E)$ be a digraph and let $\textfrak{S}_V$ be the group of permutations of $V$. Then, we define
    \[\textfrak{S}_V(X)=\{\pi\in \textfrak{S}_V\mid \textrm{ each non-trivial cycle of }\pi \textrm{ is a cycle of } X\},\]
    \[\textfrak{S}_V(X, \overline{X})=\{\pi\in \textfrak{S}_V\mid \textrm{ each cycle of }\pi \textrm{ is a cycle of } X, \textrm{ or a cycle of }\overline{X}\}.\]
\end{definition}

The following theorem shows that $U_X$ is, indeed, a symmetric function for any digraph $X$. For a permutation $\pi$, let $\mathrm{type}(\pi)$ denote the partition whose entries are the lengths of the
cycles of $\pi$. 

\begin{theorem}\cite{S}\label{pbaza} Let $X=(V, E)$ be a digraph. For any $\pi\in \textfrak{S}_V$, let $\varphi(\pi):=\sum_{\gamma}(\ell (\gamma)-1),$ where the summation runs over all cycles $\gamma$ of $\pi$ that are cycles in $X$ and $\ell(\gamma)$ denotes the length of the cycle $\gamma$. Then, \[U_X=\sum_{\pi\in\textfrak{S}_V(X, \overline{X})}(-1)^{\varphi(\pi)}p_{\mathrm{type}(\pi)}.\]
    
\end{theorem}

We say that two digraphs $X=(V,E)$ and $Y=(V',E')$ are
isomorphic if there is a bijection of vertices $f:V\rightarrow V'$
such that $(u,v)\in E$ if and only if $(f(u),f(v))\in E'$. Denote by $[X]$ the isomorphism class of a digraph $X$.

Let $\mathcal{D}=\oplus_{n\geq0}\mathcal{D}_n$ be the graded vector space over $\textbf{k}$, which is linearly spanned by the set of all isomorphism classes of digraphs, where the grading is given by the number of vertices. In \cite{GS} the space $\mathcal{D}$ is endowed with the structure of combinatorial Hopf algebra.

 For digraphs $X=(V,E)$ and $Y=(V',E')$ we define the product $X\cdot Y$ as the digraph on the disjoint union $V\sqcup V'$ with the set of directed edges \[E\cup E'\cup\{(u,v)\ |\ u\in V, v\in V'\}.\] The product of digraphs is obviously an associative, but not a commutative operation. The linear extension of the product on digraphs determines the multiplication $\mu:\mathcal{D}\otimes\mathcal{D}\rightarrow\mathcal{D}$, which turns the space $\mathcal{D}$ into a noncommutative algebra. 
 The restriction of a digraph $X=(V,E)$ on a subset $S\subseteq V$ is the
digraph $X|_S=(S,E|_S)$, where $E|_S=\{(u,v)\in E\  |\ u,v\in S\}$. The restrictions of digraphs may be used to define a comultiplication $\Delta:\mathcal{D}\rightarrow\mathcal{D}\otimes\mathcal{D}$ by

\[\Delta([X])=\sum_{W\subseteq V}[X|_{W}]\otimes[X|_{V\setminus W}],\] where $V$ is the vertex set of a digraph $X$. Evidently, this is a coassociative and a cocommutative operation. It is easy to check that $\Delta$ is an algebra morphism. The unit element is given by the class
$[\emptyset]$ of the digraph on the empty set of vertices, while the counit is given by $\epsilon([\emptyset])=1$ and $\epsilon([X])=0$ otherwise.

Let $\zeta:\mathcal{D}\rightarrow\mathbf{k}$ be the $\mathbf{k}-$algebra homomorphism given by 

\begin{equation}\label{character}
\zeta([X])=\#\{\omega\in\Sigma_V \ |\  X\mathrm{Des}(\omega)=\emptyset\}.
\end{equation}
Note that $\zeta([X])$ actually counts Hamiltonian paths of the complementary digraph $\overline{X}$ (which is an invariant of isomorphism classes of digraphs).

The fact that the space $\mathcal{D}$, together with operations defined in this way forms a CHA, sheds a new light on the Redei-Berge function.

\begin{theorem} \cite{GS} \label{univerzalni}
Let $\Psi:\mathcal{D}\rightarrow QSym$ be the universal morphism from the combinatorial Hopf algebra of digraphs to quasisymmetric functions. Then, for a digraph $X$,
\[\Psi([X])=U_X.\]
\end{theorem}

\begin{theorem} \cite{GS}\label{dual}
    For any digraph $X$, \[U_X=U_{X^{op}}.\]
\end{theorem}

\section{The Redei-Berge function of permutations and posets}

In this section, we introduce two combinatorial Hopf algebras that naturally map into the combinatorial Hopf algebra on digraphs. As a consequence, we show how the definition of the Redei-Berge function can be naturally extended to permutations and posets. This section relies on results from both \cite{S} and \cite{GS}.

\subsection{The combinatorial Hopf algebra of permutations $\mathcal{S}$}
  
Let $\textfrak{S}_n$ denote the group of permutations of the set $[n]$. An element $\sigma\in\textfrak{S}_n$ is given by a list $(\sigma(1), \sigma(2), \ldots, \sigma(n))$. By $\emptyset$, we denote the unique permutation of the empty set $\emptyset$.  We define $\mathcal{S}_n$ as a free $\textbf{k}$-module with basis $\textfrak{S}_n$. With the operations defined as follows, the vector space $\mathcal{S}=\bigoplus_{n\geq 0}\mathcal{S}_n$ has a structure of a graded connected Hopf algebra, as it is shown in \cite{ABM}.

For the unit, we take the map $u: 1\mapsto \emptyset$. The graded multiplication is given by  
$$\pi\otimes \sigma\mapsto \pi\sigma \textrm{ for }\pi \in \textfrak{S}_m, \ \sigma\in \textfrak{S}_n,$$where $\pi\sigma$ denotes the permutation $(\pi(1), \ldots, \pi(m), \sigma(1)+m, \ldots, \sigma(n)+m)$, and linearly extended on $\mathcal{S}_m\otimes\mathcal{S}_n\rightarrow\mathcal{S}_{m+n}$.

Let $M$ and $N$ be two totally ordered sets of equal size. By $\textrm{inc}_{M\rightarrow N}$, we denote the unique increasing bijection $M\rightarrow N$. If $f: A\rightarrow B$ is an injective map between two totally ordered sets, its standardization, denoted as $\textrm{st}(f)$, is defined as $\textrm{inc}_{f[A]\rightarrow[k]}\circ f\circ\textrm{inc}_{[k]\rightarrow A}$, where $k=|A|$. The graded comultiplication on $\mathcal{S}$ is then given by
$$\sigma\mapsto \sum_{S\subseteq [n]} \text{st}(\sigma|_{S})\otimes \text{st}(\sigma|_{[n]\setminus S}), \textrm{ for }\sigma \in \textfrak{S}_n,$$ and linearly extended on
$\mathcal{S}_{n}$. 
We define the counit by $\epsilon(\emptyset)=1$.

 Equipped with these operations, $\mathcal{S}$ becomes a graded connected bialgebra, and hence a graded connected Hopf algebra \cite{ABM}. It is cocommutative, but not commutative. Nevertheless, we will not follow \cite{ABM} for the definition of character. We want to define a CHA structure accordant with the structure of CHA on digraphs.

If $\sigma$ is a permutation on $[n]$, we say that a permutation $\pi\in \mathfrak{S}_n$ is $\sigma$-reversing if $\pi(i)>\pi(i+1)$, or $\sigma(\pi(i))>\sigma(\pi (i+1))$ for every $i\in[n-1]$. Finally, we take $\zeta(\sigma)$ to be the number of $\sigma$ - reversing permutations. The proof of the next lemma follows from Lemma \ref{pkarakter} and Theorem \ref{SuP} since a restriction of an algebra morphism to a subalgebra is still an algebra
morphism. Therefore, this proof could be omitted.

\begin{lemma} \label{karakter}
    If $\pi \in \textfrak{S}_m$ and $\sigma\in \textfrak{S}_n$, then $\zeta(\pi\sigma)=\zeta(\pi)\zeta(\sigma)$. Hence, $\zeta$ is a well defined character on $\mathcal{S}$.
\end{lemma}

\begin{proof}
    If $\lambda=(\lambda(1),  \ldots, \lambda(m))$ is $\pi-$reversing and $\mu=(\mu(1), \ldots, \mu(n))$ is $\sigma$-reversing, let their \textit{combination} $\tau$ be $(\mu(1)+m, \ldots, \mu(n)+m, \lambda(1), \ldots, \lambda(m))$. If $i\in [n-1]$, then the condition $\tau(i)>\tau(i+1)$, or $\pi\sigma(\tau(i))>\pi\sigma(\tau(i+1))$ is equivalent to $\mu(i)>\mu(i+1)$, or $\sigma(\mu(i))>\sigma(\mu(i+1))$, which is true since $\mu$ is a $\sigma-$reversing permutation. Similarly, for $i\in \{n+1, \ldots, n+m\}$, the condition $\tau(i)>\tau(i+1)$, or $\pi\sigma(\tau(i))>\pi\sigma(\tau(i+1))$ is actually $\lambda(i-n)>\lambda(i+1-n)$, or $\pi(\lambda(i-n))>\pi(\lambda(i+1-n))$, which is true since $\lambda$ is $\pi$-reversing. Finally, for $i=n$, we have that $\tau(i)=\mu(n)+m>m\geq \lambda(1)=\tau(i+1).$ Hence, $\tau$ is a $\pi\sigma-$reversing permutation.
    
    Conversely, if $\tau=(\tau(1), \ldots, \tau(m+n))$ is $\pi\sigma$-reversing, $\{\tau(1), \ldots, \tau(n)\}=\{m+1, \ldots, m+n\}$ and $\{\tau(n+1), \ldots, \tau(n+m)\}=\{1, \ldots, m\}$. If not, there exists $i$ such that $\tau(i)\in \{1, \ldots, m\}$ and $\tau(i+1)\in \{m+1, \ldots, m+n\}$. However, this implies that $\tau(i)\ngtr\tau(i+1)$ and that $\pi\sigma(\tau(i))\ngtr\pi\sigma(\tau(i+1))$, which is impossible since $\tau$ is $\pi\sigma-$reversing. Therefore, the condition $\pi\sigma(\tau(i))>\pi\sigma(\tau(i+1))$ is equivalent to $\sigma(\tau(i)-m)>\sigma(\tau(i+1)-m)$ for $i\in[n-1]$ and to $\pi(\tau(i-n))>\pi(\tau(i+1-n))$ for $i\in\{n+1, \ldots, n+m\}$. Hence, $\tau$ can be obtained as a combination of a $\pi$-reversing permutation $\lambda$ and a $\sigma$-reversing permutation $\mu$.
\end{proof}

\subsection{The combinatorial Hopf algebra of posets $\mathcal{P}$}

  Let $\mathcal{P}_n$ denote a free $\textbf{k}$-module with basis the set of all partial orders on the set $[n]$. The relation that induces poset $P$ will be denoted as $\leq_P$ and the corresponding strict order relation as $<_P$. For the graded vector space $\mathcal{P}=\bigoplus_{n\geq 0}\mathcal{P}_n$, it is possible to define the structure of combinatorial Hopf algebra in many ways, see \cite{ABS}. We define the operations on $\mathcal{P}$ as in \cite{ABM}. 

The unit is the map $u:\mathbf{k}\rightarrow \mathcal{P}$ given by $1\mapsto \emptyset$, where we denoted by $\emptyset$ the unique poset on the empty set. The graded multiplication $\mathcal{P}_m\otimes\mathcal{P}_n\rightarrow\mathcal{P}_{m+n}$ is the linear extension of the operation given by
 $$P\otimes Q\mapsto PQ \textrm{ for } P \textrm{ a poset on } [m] \textrm{ and $Q$ a poset on } [n]$$
where $PQ$ denotes the poset on $[m+n]$ with $\leq_{PQ}$ defined as follows: 
\begin{itemize}
\item  $x\leq_{PQ}y$ if and only if $x\leq_P y$ for all $x, y\in[m],$
\item $x+m\leq_{PQ} y+m$ if and only if $x\leq_Q y$ for all $x, y \in [n]$ and
\item $x\leq_{PQ} y$ for all $x\in [m]$, $y\in \{m+1, \ldots, m+n\}$
\end{itemize}
The poset $PQ$ defined above is known as the ordinal sum of posets $P$ and $Q$. 

For a poset $P$ on some finite set $A\subseteq\mathbb{N}$, we define the \textit{standardization} of $P$, denoted as $\text{st}(P)$, as a poset on $\{1, \ldots, |A|\}$ obtained from $P$ via the standardization of $A$. The graded comultiplication on $\mathcal{P}$ is given by 
$$P\mapsto \sum_{S\subseteq [n]} \text{st}(P|_{S})\otimes \text{st}(P|_{[n]\setminus S}),$$where $P|_{X}$ denotes the restriction of relation $\leq_P$ on $X$, and linearly extended on 
$\mathcal{P}_{n}$. We define the counit by $\epsilon(\emptyset)=1$.

With operations defined this way, $\mathcal{P}$ becomes a graded connected Hopf algebra that is cocommutative, but not commutative, see \cite{ABM}. For a poset $P$, its \textit{quasi-linear extension} is a list $(p_1, \ldots, p_n)$ of all elements of $P$ with no repetitions such that $p_{i+1}\nless p_i$ for all $i\in\{1, \ldots, n-1\}$. For the same reason as for the CHA of permutations defined above, for the character $\zeta$, we will take $\zeta(P)$ for some poset $P$ to be the number of quasi-linear extensions of $P$. Clearly, if $(r_1, \ldots, r_n)$ is a quasi-linear extension of $PQ$, then $r_1, \ldots, r_{|P|}\in P$ and $r_{|P|+1}, \ldots, r_n\in Q$. Moreover, $(r_1, \ldots, r_{|P|})$ is a quasi-linear extension of $P$ and $(r_{|P|+1}, \ldots, r_n)$ is a quasi-linear extension of Q. Therefore, $\zeta$ is multiplicative.

\subsection{The connection between CHAs $\mathcal{S}$, $\mathcal{P}$ and $\mathcal{D}$}

There is a natural connection between the three combinatorial Hopf algebras from the previous text. If $\sigma\in \textfrak{S}_n$, let $P_\sigma$ denote a poset on $[n]$ given by $i<_{P_\sigma} j$ if $i<j$ and $\sigma(i)<\sigma(j)$. Define $f_n:\mathcal{S}_n\rightarrow\mathcal{P}_n$ as a linear extension of the assignments $\sigma\mapsto P_\sigma$. This gives rise to a unique $\mathbf{k}-$linear function $f:\mathcal{S}\rightarrow \mathcal{P}$ such that $f\restriction_{\mathcal{S}_n}=f_n$ for every $n\in\mathbb{N}$. 

If $P$ is a poset on the set $[n]$, let $[D_P]$ be the isomorphism class of a digraph $D_P$ on the vertex set $[n]$ defined by the following condition: $(i, j)$ is an edge in $D_P$ if $i\leq_P j$. The linear extensions $g_n:\mathcal{P}_n\rightarrow\mathcal{D}_n$ of the assignments $P\mapsto[D_P]$ produce a unique $\mathbf{k}$-linear function $g: \mathcal{P}\rightarrow\mathcal{D}$ such that $g\restriction_{\mathcal{P}_n}=g_n$ for every $n\in\mathbb{N}$.

\begin{lemma} \label{pkarakter}
    If $\sigma\in\textfrak{S}_n$, then $\sigma-$reversing permutations are exactly the quasi-linear extensions of the poset $P_{\sigma}^{*}$, the dual of the poset $P_{\sigma}$.
\end{lemma}

\begin{proof}
Let $(p_1, \ldots, p_n)$ be a list of elements of $[n]$. This list induces a unique permutation $\pi\in\textfrak{S}_n$ by $\pi(i)=p_i$. Note that $(p_1, \ldots, p_n)$ is a quasi-linear extension of the poset $P_{\sigma}^*$ if and only if 
\[
    \begin{split}
        &(p_{i+1}\nless p_i \textrm{ in } P_{\sigma}^*)\\
        &\iff (p_i\nless p_{i+1} \textrm{ in } P_{\sigma})\\
        &\iff (p_i\geq p_{i+1} \textrm{ or }\sigma(p_i)\geq\sigma(p_{i+1}))\\
        &\iff (p_i>p_{i+1} \textrm{ or } \sigma(p_i)>\sigma(p_{i+1}))\\
        &\iff (\pi(i)>\pi(i+1) \textrm{ or }\sigma(\pi(i))>\sigma(\pi(i+1))).
    \end{split}
\]
    \end{proof}

\begin{theorem} \label{SuP}
   CHA $\mathcal{S}$ is a combinatorial Hopf subalgebra of CHA $\mathcal{P}$.
\end{theorem}

\begin{proof} The map $f: \mathcal{S}\rightarrow\mathcal{P}$ defined above is an injective morphism of graded Hopf algebras, see \cite{ABM}. Hence, we only need to check that the respective characters are compatible.

However, according to the previous lemma, $\sigma$-reversing permutations are exactly the quasi-linear extensions of $P_{\sigma}^*$, the dual poset of $P_{\sigma}$. Since a poset $P$ and its dual $P^*$ have the same number of quasi-linear extensions, we see that the characters agree.
\end{proof}

On the other hand, we cannot say that $\mathcal{P}$ is a combinatorial Hopf subalgebra of $\mathcal{D}$. This is due to the fact that for the generators of $\mathcal{P}$ we have taken all posets on the set $[n]$, while for the generators of $\mathcal{D}$ we have taken only all isomorphism classes of digraphs on the set $[n]$. Yet, if we define $\mathcal{D}$ to be generated by all digraphs on the set $[n]$, which is only a technical difference, the following theorem would yield an embedding of $\mathcal{P}$ into $\mathcal{D}$.

\begin{theorem}
    The map $g: \mathcal{P}\rightarrow\mathcal{D}$ is a morphism of combinatorial Hopf algebras.
\end{theorem}

\begin{proof}
 It is straightforward to check that $g$ respects multiplication, comultiplication, unit, counit and grading. Let $P^*$ denote the dual poset of $P$. The Hamiltonian paths of $\overline{D_P}$ are the listings $(p_1, \ldots, p_n)$ of $P$ such that $p_i\nleq_P p_{i+1}$ for all $i\in [n-1]$, which is equivalent to $p_i\nless_P p_{i+1}$ for all $i\in [n-1]$ and thus to $p_{i+1}\nless_{P^*}p_i$ for all $i\in[n-1]$. Hence, these paths are exactly the quasi-linear extensions of $P^*$. Since $P$ and $P^*$ have the same number of quasi-linear extensions, the characters are compatible.
\end{proof}

These two morphisms allow us to naturally define the Redei-Berge function for permutations and posets as the Redei-Berge function of their respective digraphs.

\begin{definition}
    Let $\sigma$ be a permutation. The Redei-Berge function $U_{\sigma}$ of $\sigma$ is \[U_{\sigma}=U_{g(f(\sigma))}.\]
Let $P$ be a poset. The Redei-Berge function $U_P$ of $P$ is \[U_P=U_{g(P)}.\] 
\end{definition}

Since $f$ and $g$ are homomorphisms of CHAs, we see immediately that, according to Theorem \ref{univerzalni}, these functions are exactly the images of $\sigma$ and $P$ by the universal homomorphism of CHAs $\mathcal{S}$ and $\mathcal{P}$ respectively to $QSym$.

\begin{theorem}
    If $\Psi: \mathcal{S}\rightarrow QSym$ is the universal morphism of combinatorial Hopf algebras, then for a permutation $\sigma$, \[\Psi(\sigma)=U_{\sigma}.\]
    If $\Psi: \mathcal{P}\rightarrow QSym$ is the universal morphism of combinatorial Hopf algebras, then for a poset $P$, \[\Psi(P)=U_{P}.\]

\end{theorem}

 In fact, we could define our morphism $g$ differently. Namely, another natural way to get a connection between posets and digraphs is to assign a digraph to a poset by taking $(u, v)$ as an edge if $v\leq_P u$. However, this would not change the Redei-Berge function of the poset, since $U_X=U_{X^{op}}$ for any digraph $X$ by Theorem \ref{dual}. Similarly, since $g(P)$ and $g(P^*)$ are two opposite digraphs, we have the following.

 \begin{lemma} \label{minmax}
     If $P$ is a poset and $P^*$ its dual poset, then 
\[U_P=U_{P^*}.\] Consequently, the number of minimal and the number of maximal elements of some poset $P$ cannot be uncovered from $U_P$. 
 \end{lemma}

 Notice that permutations and posets always induce digraphs without non-trivial directed cycles (acyclic digraphs). Also, the collection of acyclic digraphs is closed under the multiplication and comultiplication in $\mathcal{D}$ and, hence, forms a Hopf subalgebra of $\mathcal{D}$.

\begin{theorem}
    If $D$ is an acyclic digraph, then $U_D$ is $p-$positive. Consequently, for any permutation $\sigma$ and any poset $P$, $U_\sigma$ and $U_P$ are $p-$positive.
\end{theorem}
\begin{proof}
    This follows immediately from Theorem \ref{pbaza}, since if $D$ has no non-trivial cycles, $\varphi(\pi)=0$ for any permutation $\pi$. 
\end{proof}

\begin{example} \label{primer}
1) Let $L_n$ denote the natural linear order on the set $[n]$. It is obvious that $L_n=L_{n-1}\cdot L_1$, hence $U_{L_n}=\Psi(L_1^n)=\Psi(L_1)^n=U_{L_1}^n=(F_{\emptyset})^n=p_{(1, \ldots, 1)}$. Notice that we could conclude the same by using Theorem \ref{pbaza}. Since $D_{L_n}$ has no non-trivial cycles and $\overline{D_{L_n}}$ has no cycles, it follows that $\textfrak{S}_V(D_{L_n}, \overline{D_{L_n}})=\{\text{Id}\}$. In $\text{Id}$, all cycles have length 1, hence we obtain $U_{L_n}=p_{(1, \ldots, 1)}.$

2) Let $A_n$ denote the discrete order on the set $[n]$. Then, $D_{A_n}$ has no edges that are not loops. Directly from the Definition \ref{descents}, we get that $U_{A_n}=n!F_{\emptyset}$. On the other hand, since $D_{A_n}$ has only trivial cycles and $\overline{D_{A_n}}$ has all possible nontrivial cycles, $\textfrak{S}_V(D_{A_n}, \overline{D_{A_n}})=\textfrak{S}_V$, and by Theorem \ref{pbaza}, we have that $U_{A_n}=\sum_{\pi\in \textfrak{S}_V}p_{\mathrm{type}(\pi)}.$ 
\end{example}

\section{Some further properties of the Redei-Berge function}

Motivated by some properties of the chromatic symmetric function, in the sequel we will further investigate the Redei-Berge symmetric function. Subsection 4.1 is independent of the results from \cite{S} and \cite{GS}, using only the definition of the Redei-Berge function, while subsections 4.2 and 4.3 rely on both of these papers.

\subsection{Decomposition techniques for digraphs}

The Redei-Berge polynomial of a digraph, as an analogue of the chromatic polynomial of a graph, is defined in \cite{GS}. The \textit{Redei-Berge polynomial} of a digraph $X$ is the principal specialization of its Redei-Berge function: \[u_X(m)=\mathrm{ps}^1(U_X)(m)=U_X(\underbrace{1,\ldots,1}_{\text{$m$ ones}},0,0,\ldots).\]
The Redei-Berge polynomial, just like the chromatic polynomial, satisfies the deletion-contraction property, see \cite{GS}. The same property does not hold for the Redei-Berge function since the contraction of an edge gives a digraph with less vertices than the original digraph. However, there is a way to overcome this difficulty and express the Redei-Berge function of a digraph $X$ in terms of the same function of some simpler digraphs. Similar decompositions for the chromatic symmetric function can be found in \cite{DvW} and \cite{OS}. 

\begin{definition}
    A \textit{bag of sticks} is a digraph that can be written as a disjoint union of directed paths, where paths with no edges are allowed.
\end{definition}

\begin{figure}
    \centering
    \begin{tikzcd}
\bullet \arrow[d] & \bullet \arrow[d] & \bullet \arrow[d] & \bullet \\
\bullet \arrow[d] & \bullet \arrow[d] & \bullet \arrow[d] &   \\
\bullet \arrow[d] & \bullet           & \bullet           &   \\
\bullet           &             &             &  
\end{tikzcd}
    \caption{Bag of sticks}
    \label{fig:enter-label}
\end{figure}

\begin{theorem}
    Let $X=(V, E)$ be a digraph and let $F\subseteq E$ be a subset such that the subdigraph $(V, F)$ is not a bag of sticks. Then, \[\sum_{S\subseteq F}(-1)^{|S|}U_{X\setminus S}=0,\] where $X\setminus S := (V, E\setminus S)$. Consequently, \[U_X=\sum_{\substack{S\subseteq F \\ S\neq\emptyset}}(-1)^{|S|-1}U_{X\setminus S}.\]
\end{theorem}

\begin{proof} The definition of the Redei-Berge function and swapping the order of summation give us
\[\sum_{S\subseteq F}(-1)^{|S|}U_{X\setminus S}=\sum_{S\subseteq F}(-1)^{|S|}\sum_{\omega\in\Sigma_V}F_{(X\setminus S)\mathrm{Des}(\omega)}=\sum_{\omega\in\Sigma_V}\sum_{S\subseteq F}(-1)^{|S|}F_{(X\setminus S)\mathrm{Des}(\omega)}.\] We will prove that for a fixed $\omega\in \Sigma_V$, the inner sum is 0.

Let $\omega$ be a $V-$listing. Since $(V, F)$ is not a bag of sticks,  there exists an edge $e=(u, v)\in F$ such that $u, v$ do not appear consecutively in that order in $\omega$. Therefore, this edge is not detected by $X\mathrm{Des}(\omega)$. The subsets of $F$ that do not contain $e$ are in bijection with the subsets of $F$ that contain $e$ by $S\mapsto S\cup\{e\}$. Notice that $X\setminus S$ and $X\setminus (S\cup\{e\})$ induce the same descent set assigned to $\omega$, since, as we have already noted, $e$ is not detected by $X\mathrm{Des}(\omega)$. However, in $\sum_{S\subseteq F}(-1)^{|S|}F_{(X\setminus S)\mathrm{Des}(\omega)}$, terms $(-1)^{|S|}F_{(X\setminus S)\mathrm{Des}(\omega)}$ and $(-1)^{|S\cup\{e\}|}F_{(X\setminus (S\cup\{e\}))\mathrm{Des}(\omega)}$ have opposite signs, hence, they cancel each other.
\end{proof}

\begin{corollary}\label{razbijanje}
    Let $X=(V, E)$ be a digraph that is not a bag of sticks. Then, \begin{equation}\label{podskupovi}U_X=\sum_{\substack{S\subseteq E \\ S\neq\emptyset}}(-1)^{|S|-1}U_{X\setminus S}.\end{equation}
\end{corollary}

\begin{corollary}
    (k-deletion): If $e_1, e_2, \ldots, e_k$ is a list of edges that form a directed $k-$cycle in a digraph $X=(V, E)$, then
   \[\label{deletion}
U_X=\sum_{\substack{S\subseteq \{e_1, e_2, \ldots, e_k\} \\ S\neq\emptyset}}(-1)^{|S|-1}U_{X\setminus S}. 
  \]
\end{corollary}

\begin{example} \label{triple}
    (triple-deletion): If $e_1, e_2, e_3$ is a list of edges that form a triangle in a digraph $X=(V, E)$, then \[U_X=U_{X\setminus \{e_1\}}+U_{X\setminus\{ e_2\} }+U_{X\setminus\{ e_3\} }-U_{X\setminus\{ e_1, e_2\} }-U_{X\setminus\{ e_2, e_3\} }-U_{X\setminus\{ e_3, e_1\} }+U_{X\setminus\{ e_1, e_2, e_3\} }.\]
\end{example}

Previous properties can be a very powerful tool for calculating the Redei-Berge function of some particular digraph. There is no analogous feature of the chromatic function. On the other hand, the chromatic function is multiplicative on the disjoint union of two graphs, since it is an image of the universal morphism of CHA of graphs where multiplication is given by disjoint union of graphs. The same property does not hold for the Redei-Berge function.

\subsection{Properties preserved by the Redei-Berge function}

Motivated by Stanley's question from \cite{SR}, to what extent chromatic symmetric function distinguishes non-isomorphic graphs, many researchers have studied the properties of graphs that can be determined from their chromatic function, \cite{MMW}, \cite{OS} etc. In this subsection, we will try to answer some similar questions concerning the Redei-Berge function.

We do not expect that the Redei-Berge function distinguishes non-isomorphic digraphs (equivalently, that $\Psi:\mathcal{D}\rightarrow Sym$ is injective). We already noticed that $U_X=U_{X^{op}}$ and $X$ and $X^{op}$ are not necessarily isomorphic. Even for $n=2$, there is an $n-$vertex digraph $X$ that is not isomorphic to its opposite digraph : $X=(\{1, 2\}, \{(1, 1), (1, 2)\})$. Also, adding new loops does not affect the Redei-Berge function. However, we can say the following.

\begin{theorem} \label{grane}
    Two digraphs having the same Redei-Berge function have the same number of vertices and the same number of edges that are not loops.
\end{theorem}

\begin{proof}
    The first part is obvious since for digraph $X=(V, E)$, $U_X$ is homogeneous of degree $|V|$. By definition, $ U_X=\sum_{\omega\in \Sigma_V} F_{X\mathrm{Des}(\omega)}$, which can be rewritten as $\sum_{S\subseteq [n-1]}c_SF_S$ by gathering the terms that correspond to listings with the same descent set. If $(i, j)$ is an edge of $X$, there are $(n-1)!$ listings of $[n]$ in which $i$ and $j$ are consecutive. Since $|X\mathrm{Des}(\omega)|=|\{1\leq i\leq n-1\mid (\omega_i, \omega_{i+1})\in E\}|$, it follows that $(n-1)!|E|=\sum_{\omega\in\Sigma_V}|X\mathrm{Des}(\omega)|=\sum_{S\subseteq[n-1]}c_S|S|$, and therefore \[|E|=\frac{\sum_{S\subseteq[n-1]}c_S|S|}{(n-1)!}.\]
\end{proof}

\begin{corollary} \label{uporedivi}
    Two posets having the same Redei-Berge function have the same number of elements and the same number of comparable pairs.
\end{corollary}

 The \textit{length} of some path $v_1, v_2, \ldots v_k$  in a digraph is defined to be $k-1$. Similarly, the \textit{length} of some chain in a poset is the number of elements in that chain diminished by one. In fact, the previous two results can be obtained from the following for $k=0$ and $k=1$.

\begin{theorem} \label{put}
    If two digraphs with $n$ vertices have the same Redei-Berge function, they also have the same number of paths of length $k$ for any $k\in \{0, 1, \ldots, n-2\}$. Consequently, their longest paths have the same length. 
\end{theorem} 

\begin{proof}
If $v_1, v_2, \ldots, v_{k+1}$ is a path of length $k$ in $X=(V, E)$ with $|V|=n$, then the vertices $v_1, v_2, \ldots, v_{k+1}$ appear consecutively in this order in $(n-k)!$ listings of $V$. By definition, $U_X=\sum_{\omega\in \Sigma_V} F_{X\mathrm{Des}(\omega)}$, which we rewrite as $\sum_{S\subseteq [n-1]}c_SF_S$ in the same manner as in the previous proof. For $S\subseteq[n-1],$ we define $\lambda_S^k=|\{i\in [n-1]\mid i, i+1, \ldots, i+k\in S\}|.$ If the number of paths of length $k$ in $X$ is denoted by $P_X^k$, we see that \[(n-k)!P_X^k=\sum_{\omega\in\Sigma_V}\lambda_{X\mathrm{Des}(\omega)}^k=\sum_{S\subseteq[n-1]}c_S\lambda_S^k.\] Therefore, \[P_X^k=\frac{\sum_{S\subseteq[n-1]}c_S\lambda_S^k}{(n-k)!}.\]
\end{proof}

\begin{corollary} \label{posetput}
    If two posets have the same Redei-Berge function, they have the same number of chains of length $k$ for any $k\in \{0, 1, \ldots, n-2\}$. Accordingly, their longest chains have the same length.
\end{corollary}

From Corollary \ref{uporedivi}, we see that the number of incomparable pairs in some poset $P$ can be uncovered from $U_P$. This also follows directly from Theorem \ref{pbaza}, since \[[p_{(2, 1, \ldots, 1)}]U_P=\#\{(k, l)\mid k<l \textrm{ and $k$ and $l$ are incomparable in $P$}\},\]
due to the fact that $D_P$ does not contain cycles of length 2.

Now we will focus our attention on some poset $P$ that is a disjoint union of at most 4 incomparable chains $L_x, L_y, L_z, L_w$, where $x\geq y\geq z\geq w\geq 0$. From what we have said so far, we see that the following parameters can be revealed from $U_P$ : the length of the longest chain $a=x$, the number of elements $b=x+y+z+w$, the number of comparable pairs $c=\binom{x}{2}+\binom{y}{2}+\binom{z}{2}+\binom{w}{2}$ and the number of incomparable pairs $d=xy+xz+xw+yz+yw+zw$. However, since $2c=b^2-2d-b$, these four equations are dependent, and do not necessarily induce a unique solution. Nevertheless, if we set $w=0$, we get the following.

\begin{theorem}
    If two posets $P$ and $Q$ that are both disjoint unions of at most three incomparable chains have the same Redei-Berge function, then they are isomorphic.
\end{theorem}

The previous theorem can be further generalized if we exploit  Corollary \ref{posetput} in a more detailed manner.

\begin{theorem}
    If two posets $P$ and $Q$ that are both disjoint unions of incomparable chains have the same Redei-Berge function, then they are isomorphic.
\end{theorem}

\begin{proof}
    Let $(0^{r_0}, 1^{r_1}, 2^{r_2}, \ldots, k^{r_k})$ and $(0^{s_0}, 1^{s_1}, 2^{s_2}, \ldots, l^{s_l})$ be the lists of lengths of chains of $P$ and $Q$ respectively. According to Corollary \ref{posetput}, the longest chains of $P$ and $Q$ have the same length, hence $k=l$ and, furthermore, $r_k=s_k$. 
    
    If we now suppose that we have already proved that $r_i=s_i$ for every $i\in \{m+1, m+2, \ldots, k\}$, we are able to deduce that the same thing holds for $i=m$ as well. Namely, we know from Corollary \ref{posetput} that the numbers of chains of length $m$ in $P$ and $Q$ are the same. However, since there are $\binom{n+1}{m+1}$ subchains of length $m$ in any chain of length $n$ for $n\geq m$, this implies that \[r_m+\binom{m+2}{m+1}r_{m+1}+\cdots+\binom{k+1}{m+1}r_k=s_m+\binom{m+2}{m+1}s_{m+1}+\cdots+\binom{k+1}{m+1}s_k,\]
    which gives $r_m=s_m$, as desired. Hence, $r_i=s_i$ for every $i$, and therefore, $P$ and $Q$ are isomorphic.
\end{proof}
In the same manner, we can obtain the following.

\begin{theorem}
    If two digraphs that are both bags of sticks have the same Redei-Berge function, then they are isomorphic.
\end{theorem}

The importance of bags of sticks comes from the following. According to the Corollary \ref{razbijanje}, the Redei-Berge function of any digraph that is not a bag of sticks can be expressed as a linear combination of the Redei-Berge functions of appropriate subdigraphs. The same expansion can be applied again to any digraph appearing on the right side of the Equation \ref{podskupovi} that is not a bag of sticks. If we continue with this procedure, we will be able to express the Redei-Berge function of the original digraph as a linear combination of the Redei-Berge functions of its spanning subdigraphs that are bags of sticks. 

\begin{lemma} \label{span}
    $\mathrm{Span}\{U_X\mid X \text{ is a bag of sticks}\}=\mathrm{Span}\{U_X\mid X \text{ is a digraph}\}.$
\end{lemma}

Finally, there is one more thing we can say about a special class of digraphs.


\begin{definition}
    A \textit{tournament} is a loopless digraph $X=(V,E)$ such that whenever $u$ and $v$ are two distinct vertices in $V$, exactly one of $(u,v)$ and $(v,u)$ is an edge of $X$.
\end{definition}

\begin{theorem}
    If two tournaments have the same Redei-Berge function, then they have the same number of $k-$cycles for any odd $k\geq 3$.
\end{theorem}

\begin{proof}
    If $X=(V, E)$ is a tournament and $(u_1, u_2, \ldots, u_k)$ is a $k-$cycle in $X$, then $(u_k, \ldots, u_2, u_1)$ is a $k-$cycle in $\overline{X}$ and vice versa. According to Theorem \ref{pbaza}, the number of $k-$cycles in $X$ is exactly $\frac{[p_{(k, 1, \ldots, 1)}]U_X}{2}.$
\end{proof}

The previous theorem can not be expanded to the case of cycles of an even length since the terms appearing in Theorem \ref{pbaza} would cancel each other out. 

\subsection{Some new bases of $Sym$}

In this section, motivated by the result from \cite{CvW}, we provide abundance of new algebraic bases for the algebra of symmetric functions, whose generators are Redei-Berge functions. In what follows, we will suppose that we are working over a field $\mathbf{k}$ of characteristic 0 since the power sum functions form a basis of $Sym$ only in that case.

\begin{theorem} \label{baza}
    Let $(X_n)_{n\in\mathbb{N}}$ be a list of digraphs such that:

1) $X_n$ has $n$ vertices

2) at least one of $X_n$ and $\overline{X_n}$ has a Hamiltonian cycle if $n$ is odd and $n>1$

3) the numbers of Hamiltonian cycles in $X_n$ and in $\overline{X_n}$ are different if $n$ is even.

Then $(U_{X_n})_{n\in\mathbb{N}}$ is an algebraic basis of $Sym$ over $\mathbf{k}$, i.e. every symmetric function can be expressed uniquely as a polynomial over $\mathbf{k}$ in the $U_{X_n}$'s.
\end{theorem}

\begin{proof}
 Theorem \ref{pbaza} yields $[p_1]X_1=1$. For a permutation $\pi\in\textfrak{S}_V(X_n, \overline{X_n})$, $\textrm{type}(\pi)=(n)$ if and only if $\pi$ is either a Hamiltonian cycle of $X_n$, or a Hamiltonian cycle of $\overline{X_n}$. Let $h(X)$ denote the number of Hamiltonian cycles in digraph $X$. According to Theorem \ref{pbaza}, if $n$ is odd and $n>1$, $[p_n]U_{X_n}=h(X_n)+h(\overline{X_n})$.
Similarly, if $n$ is even, 
$[p_n]U_{X_n}= h(\overline{X_n})-h(X_n)$.
If the conditions 2) and 3) are fulfilled, these numbers are non-zero.

For $\lambda=(\lambda_1, \ldots, \lambda_k)\vdash n$, define $X_{\lambda}=X_{\lambda_1}\cdots X_{\lambda_k}$. The statement of this theorem is equivalent to the condition that the family $\{U_{X_{\lambda}}\}$ is a linear basis of $Sym$, since from Theorem \ref{univerzalni}
\[U_{X_\lambda}=U_{X_{\lambda_1}}\cdots U_{X_{\lambda_k}}.\]
From this formula and Theorem \ref{pbaza} it follows that in the expansion of $U_{X_{\lambda}}$ in the $p-$basis, there can only appear $p_{\mu}$ for $\mu\leq\lambda$. Here, $\mu\leq\lambda$ denotes the usual dominance order on partitions, given by $\mu_1+\cdots+\mu_j\leq\lambda_1+\cdots+\lambda_j$ for any $j$.

Therefore, $(U_{X_{\lambda}})_{\lambda\vdash n}$ induces a triangular transition matrix to the power sum basis, with non-zero entries on the main diagonal. Hence, $\{U_{X_{\lambda}}\}$ is a linear basis of $Sym$ as well.
\end{proof} 
It is worth mentioning that the second half of this proof could also be done using Exercise 2.5.24 from \cite{GR}. The conditions given in the previous theorem are not too strong. For example, we could take a discrete digraph on $n$ vertices as $X_n$. 

\begin{theorem}
    There are infinitely many lists $(X_n)_{n\in\mathbb{N}}$ of loopless digraphs such that $(U_{X_n})_{n\in\mathbb{N}}$ is an algebraic basis of $Sym$. 
\end{theorem}

\begin{proof}
    For $n\in\mathbb{N}$, where $n\neq 1$, we are going to construct $n-1$ non-isomorphic digraphs with $n$ vertices such that each of these digraphs satisfies the conditions of the previous theorem. Moreover, their Redei-Berge functions will also be different. 

    For $k\in \{0, 1, \ldots, n-2\}$, let $X_{n, k}$ denote the digraph with the vertex set $[n]$ and with edges $(1, 1+i)$ for $1\leq i\leq k$, see Figure 2. This digraph, obviously, does not contain a Hamiltonian cycle, but its complement $\overline{X_{n, k}}$ has at least one such cycle - $(n, n-1, \ldots, 2, 1, n)$. Hence, $X_{n, k}$ satisfies the conditions of the previous theorem for any $k\in \{0, 1, \ldots, n-2\}$. Moreover, since $X_{n, k}$ has exactly $k$ edges, by Theorem \ref{grane}, $X_{n, k}$ and $X_{n, l}$ for $k\neq l$ have different Redei-Berge functions.
\end{proof}

\begin{theorem}
    The set $\{U_X\mid X \text{ is a bag of sticks}\}$ is a basis of $Sym$.
\end{theorem}

\begin{proof}
Let us denote by $Sym^n$ the linear subspace of $Sym$ that consists of homogeneous symmetric functions of degree $n$.
According to Lemma \ref{span} and Theorem \ref{baza}, we have that $Sym=\mathrm{Span}\{U_X\mid X \text{ is a digraph}\}=\mathrm{Span}\{U_X\mid X \text{ is a bag of sticks}\}$. Hence, we only need to prove that the Redei-Berge functions of bags of sticks are linearly independent. However, this follows immediately from the fact that the number of non-isomorphic bags of sticks with $n$ vertices is $P(n)$, the number of partitions of $n$, which is exactly the dimension of $Sym^n$.
\end{proof}

\begin{example}
    Theorem \ref{pbaza} easily gives us expansions of the Redei-Berge function of some simple bags of sticks in the power sum basis:
    \[\begin{split}
        &U_{\bullet}=p_1;\\
        &U_{\bullet\hspace{2mm} \bullet}=p_{1, 1}+p_2;\\
        &U_{\bullet\rightarrow\bullet}=p_{1, 1};\\
        &U_{\bullet\hspace{2mm}\bullet\hspace{2mm}\bullet}=p_{1, 1, 1}+3p_{2, 1}+2p_{3};\\
        &U_{\bullet\rightarrow\bullet\hspace{2mm}\bullet}=p_{1, 1,1}+2p_{2, 1}+p_3;\\
        &U_{\bullet\rightarrow\bullet\rightarrow\bullet}=p_{1, 1, 1}+p_{2, 1}+p_3.
    \end{split}\]
\end{example}

\begin{example}
    The Redei-Berge function of a directed 3-cycle (triangle) $X=(\{1, 2, 3\}, \{(1,2), (2, 3), (3, 1)\})$ could easily be calculated in the same manner - using Theorem \ref{pbaza}. However, triple-deletion from Example \ref{triple} gives another way to do so:\[U_X=3U_{\bullet\rightarrow\bullet\rightarrow\bullet}-3U_{\bullet\rightarrow\bullet\hspace{2mm}\bullet}+U_{\bullet\hspace{2mm}\bullet\hspace{2mm}\bullet}=p_{1, 1, 1}+2p_3.\]
\end{example}

\begin{figure} \label{XNK}
    \centering
    \begin{tikzcd}
\bullet \arrow[r] \arrow[rr, bend left] \arrow[rrr, bend right] & \bullet & \bullet & \bullet & \bullet
\end{tikzcd}
    \caption{$X_{5, 3}$}
    \label{fig:enter-label}
\end{figure}

\section{Further avenues}

The Redei-Berge function is a symmetric function associated to digraphs, introduced in the 1990s. Unlike its peer, the chromatic symmetric function, it did not attract attention of many researchers until recently. In this paper, we tried to derive some of its properties that, to a certain extent, resemble properties of the chromatic symmetric function.

A natural avenue to pursue is to give an expansion of the Redei-Berge function in other bases of the space of symmetric functions and to examine its positivity in these bases. For now, it is not quite clear how this could be achieved. One may use formulas given in Theorems \ref{pbaza} and \ref{razbijanje} to obtain such expansions. Though, the path of finding combinatorial interpretations of coefficients in other bases might be completely independent of these formulas. This would enable us to explicitly calculate the Redei-Berge function and the corresponding polynomial of a larger class of digraphs. We could also try to get an insight into the relationship between the standard bases of $Sym$ and the new bases described in the previous section. Finally, we can seek some further properties of digraphs that are preserved by the Redei-Berge function.

As we have already mentioned, permutations and posets necessarily induce acyclic digraphs. This fact could imply some additional nice properties of their Redei-Berge function. In \cite{SR}, Stanley gives combinatorial connection between the coefficients of the chromatic function of a graph in its expansion in elementary basis and the number of acyclic orientations of that graph, with a certain number of \textit{sinks}. In our terminology, the sinks represent the maximal elements of the poset, and their counterparts represent its minimal elements. However, we have seen in Lemma \ref{minmax} that the number $m_P$ of minimal and the number $M_P$ of maximal elements of poset $P$ cannot be revealed from $U_P$. Nevertheless, we were not able to find a counterexample dealing with the value $m_P+M_P$.

\begin{conjecture}
If two posets $P$ and $Q$ have the same Redei-Berge function, then $m_P+M_P=m_Q+M_Q$.
\end{conjecture}




\bibliographystyle{model1a-num-names}
\bibliography{<your-bib-database>}




\section{Declarations and statements}

The second researcher was supported by the Science Fund of the Republic of Serbia, grant
no. 7749891. Graphical Languages - GWORDS. The authors have no relevant financial or
 non-financial interests to disclose. All authors contributed to the study concep
tion and design. All authors read and approved the final manuscript

\section{Data availability}

No data was used for the research described in the article.




\end{document}